%% file: Wills5_5.tex
\documentclass[12pt, reqno]{amsart}
\usepackage[utf8]{inputenc}
\usepackage{amsmath}
\usepackage{amsthm}
\usepackage{amssymb}
\usepackage[numbers, square]{natbib}
\theoremstyle{plain}
\newtheorem{theorem}{Theorem}[section]
\newtheorem{lemma}[theorem]{Lemma}
\newtheorem{corollary}[theorem]{Corollary}
\newtheorem{proposition}[theorem]{Proposition}

\newtheorem*{theoremA}{Theorem A}
\newtheorem*{theoremB}{Theorem B}
\newtheorem*{theoremC}{Theorem C}

\theoremstyle{definition}
\newtheorem{definition}[theorem]{Definition}
\newtheorem{example}[theorem]{Example}

\theoremstyle{remark}
\newtheorem{remark}[theorem]{Remark}

\newcommand{\eps}{\varepsilon}
\newcommand{\la}{\langle}
\newcommand{\ra}{\rangle}
\newcommand{\conv}{\operatorname{conv}}
\newcommand{\diam}{\mathop{\mathrm{diam}}\nolimits}
\newcommand{\dist}{\mathop{\mathrm{dist}}\nolimits}
\newcommand{\Vol}{\mathop{\mathrm{Vol}}\nolimits}

\usepackage[top=3cm,bottom=3cm, left=3cm, right=3cm]{geometry}

\begin{document}

\author[M.~Dospolova]{Mariia Dospolova}
\address{Mariia Dospolova, 
St.~Petersburg Department of Steklov Institute of Mathematics, Fontanka~27, 191011 St.~Petersburg, Russia; and Faculty of Mathematics, Bielefeld University, Bielefeld, Germany}
\email{dospolova.maria@yandex.ru}

\author[M.~Germanskov]{Mikhail Germanskov}
\address{Mikhail Germanskov, St.~Petersburg Department of Steklov Institute of Mathematics, Fontanka~27, 191011 St.~Petersburg, Russia}
\email{mgermanskov@gmail.com}

\author[D.~Zaporozhets]{Dmitry Zaporozhets}
\address{Dmitry Zaporozhets, St.~Petersburg Department of Steklov Institute of Mathematics, Fontanka~27, 191011 St.~Petersburg, Russia}
\email{zap1979@gmail.com}

\title[Regularity of Gaussian processes]{Regularity of Gaussian processes, metric entropy and the Wills functional}
\keywords{Wills functional, metric entropy, Gaussian processes, sample path regularity, intrinsic volumes, Steiner entire function}
\subjclass[2020]{52A05, 30D15, 60G15, 60G17}
\thanks{This work was supported by the Ministry of Science and Higher Education of the Russian Federation (agreement 075-15-2025-344 dated 29/04/2025 for Saint Petersburg Leonhard Euler International Mathematical Institute).}

\begin{abstract}
We establish exact relations between three quantities: the growth of the Gaussian Wills functional under dilations, the small-scale metric entropy of the index set, and the critical H\"older exponent of the associated isonormal Gaussian process. For every
compact Gaussian-bounded set $T$, the
Wills growth exponent $\rho(T)$, the entropy exponent $\eta(T)$, and the
critical H\"older exponent $\alpha_H(T)$ satisfy
\[
\eta(T)=\frac{2\rho(T)}{2-\rho(T)},
\qquad
\alpha_H(T)
=1-\frac{\eta(T)}{2}
=\frac{2-2\rho(T)}{2-\rho(T)}.
\]
Thus any one of these parameters determines the other two. Our main tools
are a Wills-functional analogue of Sudakov's minoration and a complementary
upper bound for $W(rT)$ in terms of local Gaussian widths.
\end{abstract}

\maketitle

\input{introduction}
\input{preliminaries}
\input{main_results}
\input{appendix}

\input{references}
\end{document}

%% file: introduction.tex
\section{Introduction}

\subsection{Regularity of the isonormal process}

Let $H$ be a separable real Hilbert space. The \emph{isonormal process} on a subset $T \subset H$ is the centered Gaussian process $\left(X_t\right)_{t \in T}$ whose covariance structure agrees with the inner product of $H$:
\[\mathbb{E} X_t X_s = \left\langle t, s \right\rangle,\]
for all $t, s \in T$. Its canonical metric is the Hilbert distance:
\[
\left(\mathbb{E}|X_t-X_s|^2\right)^{1/2}=\|t-s\|.
\]
Many classical Gaussian processes arise as restrictions of the isonormal process to a suitable subset of a Hilbert space. For instance, Brownian motion on $[0,1]$ is the restriction of the isonormal process on $L^2[0,1]$ to the curve $\left\{\mathbf{1}_{[0,t]} : t \in [0,1]\right\}$. This curve is known as the Wiener spiral; we return to it in Example~\ref{ex:wiener-spiral}.

In \cite{Dud67}, Dudley proposed to study subsets of a Hilbert space through the regularity of the sample paths of the isonormal process indexed by them. Two regularity classes are central to this program. A set $T$ is called \emph{Gaussian bounded} (a \emph{GB-set}) if the isonormal process admits a modification whose sample paths are almost surely uniformly bounded on $T$. It is called \emph{Gaussian continuous} (a \emph{GC-set}) if the process admits a modification with almost surely continuous sample paths. Necessary and sufficient conditions for these properties were investigated by Dudley, Sudakov, Fernique, and Talagrand, among others \cite{Dud67, Sud71, Sud76, Fer75, Tal87}. Sudakov announced in \cite{Sud71}, with detailed proofs in \cite{Sud76}, the following criterion: $T$ is a GB-set if and only if the quantity
\[\mathbb{E} \sup_{t \in T} X_t := \sup_{T_0 \subset T \text{ finite}} \mathbb{E} \max_{t \in T_0} X_t,\]
called the \emph{expected supremum} of the process, is finite.

A fundamental problem, going back to \cite{Dud67}, is to bound expected suprema from above and below in terms of the intrinsic geometry of the index set;
we refer to Talagrand's monograph \cite{Tal21} for a comprehensive account. For a compact metric space $(T,d)$, let $N(T,\eps)$ denote the covering number at scale $\eps$, i.e., the smallest number of $\eps$-balls with centers in $T$ covering $T$,
and let $M(T,\eps)$ denote the packing number, i.e., the largest cardinality of an $\eps$-separated subset. The quantity $\log N(T,\eps)$ is called the metric entropy of $T$ at scale $\eps$.
The two classical estimates are Dudley's entropy upper bound \cite{Dud67} and Sudakov's minoration \cite{Sud71}.
By Dudley's theorem, finiteness of $\int_0^{\diam T}\sqrt{\log N(T,\eps)}\,d\eps$ guarantees a uniformly continuous modification together with an explicit upper modulus of continuity. A complementary lower bound is provided by Sudakov's minoration, which bounds the expected supremum from below by $c\eps\sqrt{\log M(T,\eps)}$ with a universal constant $c$.

However, entropy estimates alone characterize neither Gaussian boundedness nor Gaussian continuity; a complete description is given by Talagrand’s majorizing measure theorem \cite{Tal87, Tal21}.
For convex index sets, a complementary approach rooted in the theory of mixed volumes \cite{Sch14} connects expected suprema with convex geometry \cite{Chev76, Tsi85, Vit96, Vit01}.

\subsection{Intrinsic volumes and the Wills functional}

The intrinsic volumes $V_j(K)$ of a convex body generalize volume, surface area, and mean width. Their precise definitions and all properties needed in this paper are recalled in Subsection~\ref{subsec:intrinsic-volumes}.
In finite dimensions they are the coefficients of the Steiner polynomial, and
their independence of the ambient dimension permits an extension to compact
convex subsets of a Hilbert space. In particular, the intrinsic volumes are
defined for compact convex GB-sets in infinite dimensions as well. For a compact convex set
$K$, the first intrinsic volume satisfies Sudakov's formula
\[
V_1(K)=\sqrt{2\pi}\,\mathbb{E}\sup_{t\in K}X_t.
\]
Thus $K$ is GB if and only if $V_1(K)<\infty$, in which case all its intrinsic
volumes are finite.

If $K$ is a compact convex GB-set, one can define its \emph{Wills functional}
as the sum of all intrinsic volumes:
\[W(K) := \sum_{j=0}^\infty V_j(K).\]
The Wills functional was introduced by Wills \cite{Wil73} in connection with
lattice-point enumeration and was subsequently studied by Hadwiger
\cite{Had75}. Vitale \cite{Vit96} adapted it to Gaussian processes, while
further geometric inequalities and a generalized version of the functional
were developed in \cite{AHY21}.

This functional also admits the integral representation: for any
compact convex set $K \subset \mathbb{R}^d$,

\begin{equation}\label{eq:integral-repr-wills}
W(K) = \int_{\mathbb{R}^d} \exp\left(-\pi \dist^2(x, K)\right)\,dx,
\end{equation}
where $\dist(x, K) = \inf_{y \in K} \|x - y\|$; see \cite{Vit96, Bog98}. This
representation expresses $W(K)$ as a Gaussian-weighted tube volume.
Tsirelson's infinite-dimensional Steiner formula \cite[Theorem 5]{Tsi85}
provides a Gaussian representation of $W(K)$ and its dilations: for every
compact convex GB-set $K\subset H$ and every $r\geqslant0$,
\begin{equation}\label{eq:tsirelson-steiner}
\mathbb{E}\exp\Bigl(\sup_{t\in K}
\bigl[\sqrt{2\pi}\,r X_t-\pi r^2\|t\|^2\bigr]\Bigr)
=\sum_{j=0}^{\infty}r^j V_j(K).
\end{equation}

The sum defining $W(K)$ is finite for every compact convex GB-set. Indeed,
Chevet \cite{Chev76} and McMullen \cite{McM91} proved that
\[V_j(K) \leqslant \frac{\left(V_1(K)\right)^j}{j!},\]
for all $j \geqslant 0$.
This bound follows from the ultra-log-concavity of
the intrinsic volumes (see Proposition~\ref{prop:intrinsic-volume-inequalities}), which in turn is a consequence of the
Alexandrov--Fenchel inequality (see \cite[Section 4.2]{BZ88} or
\cite[Section 7]{Sch14}).

This yields the exponential bound
$W(K) \leqslant \exp\left(V_1(K)\right) < \infty$
for the Wills functional.

\subsection{The Steiner entire function and the order $\rho$}

The right-hand side of~\eqref{eq:tsirelson-steiner} can be viewed as a power
series in $r$.
For a convex body $K \subset \mathbb{R}^d$, Hern\'andez Cifre and Yepes
Nicol\'as \cite{HCYN13} studied the \emph{Wills polynomial}
$\sum_{k=0}^d V_k(K) z^k$ and its roots, continuing a line of research on the
roots of Steiner polynomials \cite{HHC08}. For a compact convex GB-set
$K \subset H$, the corresponding generating series is the \emph{Steiner entire
function}
\[f_K(z) := \sum_{j=0}^{\infty} V_j(K) z^j,\]
introduced and studied in \cite{DGZ25}. The Chevet--McMullen inequality
guarantees that the series converges on the whole complex plane and satisfies
$\left|f_K(z)\right| \leqslant \exp\left(V_1(K) |z|\right)$. Its exponential order
\[\rho(K) = \limsup_{r \to \infty} \frac{\log \log M_{f_K}(r)}{\log r}, \qquad M_{f_K}(r) := \max_{|z| = r} \left|f_K(z)\right|,\]
always satisfies $0 \leqslant \rho(K) \leqslant 1$, and the whole range
$[0, 1]$ is attained \cite[Theorem 3.1]{DGZ25}. We adopt the usual convention
that a constant entire function has order zero. (We refer to \cite{Lev96} for
the basic notions of the theory of entire functions; the facts we use are
collected in Appendix~\ref{app:entire}.)
The order $\rho(K)$ reflects regularity properties of the isonormal process on
$K$. Building on Vitale's characterization of Gaussian continuity in terms of
the ratios of consecutive intrinsic volumes \cite{Vit01, Vit09}, it was shown
in \cite{DGZ25} that $K$ is a GC-set if and only if either $\rho(K) < 1$, or
$\rho(K) = 1$ and $f_K$ has type $0$. In the boundary case $\rho(K) = 1$, the
type of $f_K$ equals the It\^o--Nisio oscillation of the process on $K$ in
Vitale's normalization \cite{DGZ25,IN68,Vit01}. As an application, \cite{DGZ25} gave a counterexample to a conjecture of Gao and Vitale \cite{GV01} on the possible decay rates of the intrinsic volume ratios.

\subsection{Main results}

It is natural to ask whether the order $\rho(K)$ determines not only whether the sample paths are continuous, but also how regular they are. We show that it does, by relating three asymptotic characteristics of an index set: the entropy exponent, the critical H\"older exponent of the sample paths, and the Wills growth exponent.

We now work with compact index sets that need not be convex, using $T$ for
general nonempty compact subsets of $H$ and reserving $K$ for convex compact
sets. To avoid degenerate cases, we assume throughout that every index set contains at least two points. By analogy with the order of an entire function, we define the Wills growth exponent by
\[
\rho(T):=\limsup_{r\to\infty}\frac{\log\log W(rT)}{\log r},
\]
where $W$ is the Gaussian Wills functional, the probabilistic extension of
the left-hand side of~\eqref{eq:tsirelson-steiner} to general index sets.
See~\eqref{eq:wills-gaussian} for the definition of $W$ and
Definition~\ref{def:wills-growth} for the precise conventions in the
definition of $\rho(T)$. For compact convex GB-sets, $\rho(K)$ coincides
with the order of the Steiner entire function.

Following Dudley~\cite{Dud67}, for a compact set $T$
define its \emph{entropy exponent} by
\[
\eta(T):=\limsup_{\eps\to0}
\frac{\log\log N(T,\eps)}{\log(1/\eps)}.
\]
We also define its \emph{critical H\"older exponent}, with respect to the
Hilbert metric, by
\[
\alpha_H(T)
:=
\sup\Bigl(
\{0\}\cup
\{\beta>0:X\text{ has a }\beta\text{-H\"older modification on }T\}
\Bigr).
\]

Our main result gives explicit relations connecting these three quantities.

\begin{theoremA}[Theorem~\ref{thm:asymptotic-exponents}]
Let $T\subset H$ be a compact GB-set. Then
\[
\eta(T)=\frac{2\rho(T)}{2-\rho(T)},
\qquad
\rho(T)=\frac{2\eta(T)}{2+\eta(T)},
\]
and
\[
\alpha_H(T)
=1-\frac{\eta(T)}{2}
=\frac{2-2\rho(T)}{2-\rho(T)}.
\]
\end{theoremA}

Thus the Wills growth exponent, the entropy exponent, and the critical H\"older exponent determine one another, even when the index set is not convex. The idea of considering the Wills functional for sets that are not necessarily convex first appeared in \cite{Mou25}, and it will be essential for Theorem~B below. The nonconvex case is especially interesting for processes typically parametrized by a real parameter $t\in\mathbb{R}$. The index sets of such processes are represented by curves in a Hilbert space. The Gaussian form of the Wills functional allows us to estimate their entropy growth and moduli of continuity.

In the convex case, since the order of
$f_K$ is determined by the asymptotic behavior of its Taylor coefficients,
Theorem~A relates the decay of the intrinsic volumes to metric entropy and to
the sample-path regularity of the isonormal process. We will use this relation between the order and the
Taylor coefficients in Example~\ref{ex:brownian-motion-body}, concerning the convex
hull of the Wiener spiral, and in the discussion that follows; see
Remark~\ref{rem:tsirelson-milman-pisier}.

The proof of Theorem~A rests on two quantitative inequalities, stated below as
Theorems B and C; both may be of independent interest. Theorem B is a
Wills-functional analogue of Sudakov's minoration: it bounds the metric entropy
of $T$ in terms of the growth of $W(rT)$.

\begin{theoremB}[Theorem~\ref{thm:entropy-from-wills}]
Let $T\subset H$ be a nonempty compact set. Then there exists a universal
constant $C>0$ such that, for every $\eps>0$,
\[
\sqrt{\log N(T,\eps)}
\leqslant
C\inf_{r>0}\left(
\frac{\log W(rT)}{r\eps}+r\eps
\right).
\]
\end{theoremB}

The complementary inequality, Theorem~C, bounds $\log W(rT)$ from above in
terms of the local Gaussian width function
\[
w_T^{\mathrm{loc}}(\eps):=
\sup_{t_0\in T}
\mathbb{E}\sup_{t\in T\cap B(t_0,\eps)}
\left(X_t-X_{t_0}\right).
\]
Here $B(t_0,\eps)$ is the closed Hilbert ball of radius $\eps$ centered at
$t_0$.

\begin{theoremC}[Theorem~\ref{thm:wills-from-local-width}]
There exists a universal constant $C>0$ such that, for every nonempty compact
$T\subset H$ and every $r>0$,
\[
\log W(rT)
\leqslant
\inf_{\eps>0}\left\{
C\int_{\eps}^{\infty}
\frac{w_T^{\mathrm{loc}}(t)^2}{t^3}\,dt
+\sqrt{2\pi}\,r w_T^{\mathrm{loc}}(\eps)
\right\}.
\]
\end{theoremC}

Substituting power-law majorants into Theorems B and C allows the asymptotic exponents
to be compared with one another; these inequalities lead to
the exact identities of Theorem A.

\begin{remark}
Our main result relates the critical power-law exponents for metric entropy,
the Wills functional, and the modulus of continuity.
Theorems~B and~C also allow one to compare more general majorants, such as
$t^a\lvert\log t\rvert^b$, and thus obtain estimates that retain logarithmic
factors; see Example~\ref{ex:wiener-spiral} for an illustration.
We state the corresponding equivalence in the power-logarithmic scale
without proof; a detailed treatment of more general classes of majorants
will be given in future work.
For a nonempty compact set $T\subset H$, $\rho_1\in(0,1)$, and
$\rho_2\in\mathbb{R}$, one has
\[
\log W(rT)=O\bigl(r^{\rho_1}(\log r)^{\rho_2}\bigr)\quad (r\to\infty)
\]
if and only if
\[
\log N(T,\eps)=O\bigl(\eps^{-\eta_1}(\log(1/\eps))^{\eta_2}\bigr)\quad (\eps\to0),
\]
where
\[
\eta_1=\frac{2\rho_1}{2-\rho_1},
\qquad
\eta_2=\frac{2\rho_2}{2-\rho_1}.
\]

Power-logarithmic rates arise naturally in the study of small ball
probabilities for Gaussian processes. For a compact Gaussian-continuous
index set $T\subset H$, the entropy duality theorem of Artstein, Milman,
and Szarek~\cite{AMS04} compares the metric entropy of the closed absolutely
convex hull $K=\overline{\conv}(T\cup(-T))$ in the Hilbert norm with that of
the dispersion ellipsoid (the unit ball of the reproducing kernel Hilbert
space) in the supremum norm, up to universal constants and rescaling of
the covering radius. The entropy approach of Kuelbs and Li~\cite{KL93},
further developed by Li and Linde~\cite{LL99}, relates the latter entropy
to small ball probabilities; see also Lifshits~\cite[Subsection~11.5]{Lif12} for an
exposition. Combined with the equivalence above applied to $K$, these
results relate small ball estimates in the power-logarithmic regime to
the growth of $W(rK)$, and hence to the asymptotic behavior of the
intrinsic-volume sequence of $K$.

During the final preparation of this manuscript, we became aware of
Mourtada's preprint~\cite{Mou26}, which studies the connection between
Gaussian small ball probabilities, metric entropy, and intrinsic volumes.
In view of this overlap, we restrict ourselves to noting this connection
here and do not pursue it in detail.  
\end{remark}

Dudley~\cite{Dud67} introduced the \emph{volume exponent}
$\operatorname{EV}(K)$ in terms of maximal volumes of finite-dimensional
orthogonal projections. For centrally symmetric convex sets, Milman and
Pisier proved Dudley's conjecture \cite[Conjecture~5.4]{Dud67} on Gaussian
continuity under a decay assumption on these projection volumes. Their theorem also gives the
identity \cite[Theorem~1]{MP87}
\[
\operatorname{EV}(K)+\frac12=-\frac1{\eta(K)}.
\]
Combined with Theorem~A, this shows that for centrally symmetric compact
convex GB-sets with $0<\rho(K)<1$, the limsup in the definition of
$\operatorname{EV}(K)$ remains unchanged when maximal projection volumes are
replaced by the intrinsic volumes: both equal $-1/\rho(K)$.
The corresponding Gaussian-continuity
statement formulated in terms of the intrinsic-volume exponent is precisely
the result of \cite{DGZ25} mentioned at the end of the preceding subsection.
While the maximal projection volumes introduced by Dudley already capture substantial information about sample-path regularity, the intrinsic volumes are arguably a more natural object, as their asymptotics distinguish between the classes of GB- and GC-sets see e.g., \cite{Vit01}.

\subsection{Relation to earlier work and organization of the paper}

Following Dudley~\cite{Dud67}, one can measure the size of a compact index set
and the regularity of the corresponding Gaussian process through asymptotic
exponents and study the relations among them. Our contribution is to bring the
dilation growth of the Gaussian Wills functional into this framework: in the
convex case this is the decay of the intrinsic volumes, equivalently the order
of the Steiner entire function.

The paper is organized as follows. In Section~\ref{sec:prelim} we collect the
preliminaries on covering and packing numbers, the entropy exponent, moduli of
continuity, local Gaussian widths, intrinsic volumes, and the Gaussian Wills
functional. Section~\ref{sec:main} contains the proofs of Theorems~A, B, and~C and
the application to projection volumes via the Milman--Pisier theorem. The proofs
of Corollaries~\ref{cor:entropy-upper-bound}, \ref{cor:holder},
and~\ref{cor:order-from-local-width} are given in Appendix~\ref{app:proofs},
and the proof of Theorem~\ref{thm:width-holder-equivalence} is given in
Appendix~\ref{app:width-regularity}. Appendix~\ref{app:examples} presents
several examples, and Appendix~\ref{app:entire} recalls the formulas for the
order and type of an entire function in terms of its Taylor coefficients.

\subsection{Acknowledgements and AI statement}

The authors thank Mikhail Lifshits for his careful reading of the paper and valuable comments. The authors also thank Nikita and Egor Dobronravov, who independently proposed a conjecture about the relationship between the order of the Steiner entire function and the critical H\"older exponent.

All the results in this manuscript were obtained during the spring of 2026. AI tools were only used for assistance in writing and editing the text, searching for references, and proofreading for typos.

%% file: preliminaries.tex
\section{Preliminaries}\label{sec:prelim}

Throughout, $H$ is a separable real Hilbert space, and all metric quantities on
subsets of $H$ are computed with respect to the Hilbert norm.
Further, we use separable modifications of the isonormal process on their respective index sets.  For the isonormal process such modifications exist because its canonical metric is the Hilbert distance and the index sets are separable (see~\cite[Section~4, Proposition~1]{Lif95}). In particular, the existence of a separable modification implies that all suprema considered here are measurable.
When a continuous
modification is obtained, we denote it again by $X$.

\subsection{Covering and packing numbers}
For a compact metric space $(T,d)$ and $\eps > 0$, let $N(T,\eps)$ denote the
least number of closed balls of radius $\eps$, with centers in $T$, whose union
contains $T$, and let $M(T,\eps)$ denote the largest cardinality of an
$\eps$-separated subset of $T$, that is, a subset $A \subset T$ such that
$d(x,y) \geqslant \eps$ whenever $x,y \in A$ and $x \neq y$. Compactness
ensures that both quantities are finite. We use repeatedly the comparison
\begin{equation}\label{eq:packing-covering}
N(T, \eps) \leqslant M(T, \eps),
\end{equation}
and the scaling relations
\begin{equation}\label{eq:scaling}
M(rT, r\eps) = M(T, \eps), \qquad N(rT, r\eps) = N(T, \eps),
\end{equation}
for every nonempty compact $T\subset H$ and $r,\eps>0$. Throughout the paper,
$B(t_0,\eps):=\{t\in H:\|t-t_0\|\leqslant\eps\}$ denotes the closed ball of
radius $\eps$ centered at $t_0$.

Following Dudley~\cite{Dud67}, for a
compact set $T$ we define the \emph{entropy exponent}
\begin{equation}\label{eq:entropy-exponent}
\eta(T):=\limsup_{\eps\to0}
\frac{\log\log N(T,\eps)}{\log(1/\eps)}.
\end{equation}
We use Sudakov's minoration in the following standard form.

\begin{proposition}[Sudakov minoration~{\cite{Sud71}}]\label{prop:sudakov-minoration}
  There is a universal constant $c>0$ such that, for every nonempty compact
  $T\subset H$ and every $\eps>0$,
  \[\mathbb{E}\sup_{t\in T}X_t \geqslant c\eps\sqrt{\log M(T,\eps)}
  \geqslant c\eps\sqrt{\log N(T,\eps)}.\]
\end{proposition}
For a compact GB-set, Proposition~\ref{prop:sudakov-minoration} implies
$\eta(T)\leqslant2$.

\subsection{Uniform sample moduli and Dudley's theorem}
Following the terminology of~\cite[\S15]{Lif95}, let $X=(X_t)_{t\in T}$ be a
random function on a compact metric space $(T,d)$. Writing $X(\omega)$ for the
sample path $t\mapsto X_t(\omega)$, set
\[
\theta(X(\omega),\eps)
:=
\sup_{\substack{s,t\in T\\d(s,t)\leqslant\eps}}
|X_s(\omega)-X_t(\omega)|,
\qquad \eps>0.
\]
For each fixed $\eps$, this is a random variable under the separability
convention made above. Let
$\Theta\colon[0,\infty)\to[0,\infty)$ be nondecreasing and subadditive in a neighborhood of $0$, with
$\Theta(0)=0$, $\Theta(\eps)>0$ for $\eps>0$, and
$\Theta(\eps)\to0$ as $\eps\downarrow0$. We call the function $\Theta$ a \emph{uniform
sample modulus} for $X$ if $X$ has a modification $\widetilde X$ such that
\[
\limsup_{\eps\downarrow0}
\frac{\theta(\widetilde X(\omega),\eps)}{\Theta(\eps)}
<\infty
\]
almost surely. In other words, $\Theta(\eps)$ bounds the modulus of continuity
of almost every sample path for small $\eps$, up to a random multiplicative
constant.

For the isonormal process, the canonical metric is the Hilbert distance. If
\[
\lim_{\eps\downarrow0}\frac{\eps}{\Theta(\eps)}=0,
\]
then Marcus and Rosen~\cite[Lemma~7.1.1]{MR06} show that the upper limit above
is almost surely equal to a deterministic constant
$C_\Theta\in[0,\infty]$. Hence $C_\Theta<\infty$ whenever $\Theta$ is a
uniform sample modulus. This applies, in particular, when
$\Theta(\eps)=\eps^\alpha$ for $0<\alpha<1$.

The Dudley entropy integral provides a canonical uniform sample modulus.

\begin{theorem}[Dudley~{\cite{Dud73}; see also~\cite[\S15, Theorem~1]{Lif95}}]
\label{thm:dudley}
Let $T\subset H$ be compact and suppose that
\[
\mathcal{D}(T,\delta)
:=
\int_0^\delta\sqrt{\log N(T,\eps)}\,d\eps
\]
is finite for some (and hence every) $\delta>0$. Then the isonormal process
on $T$ has a modification $\widetilde X$ with almost surely uniformly
continuous sample paths, and $\mathcal{D}(T,\cdot)$ is a uniform sample
modulus for $\widetilde X$. More precisely,
\[
\limsup_{\delta\downarrow0}
\frac{\theta(\widetilde X(\omega),\delta)}
{\mathcal{D}(T,\delta)}
\leqslant 36
\]
almost surely.
\end{theorem}

\subsection{The H\"older modulus of continuity}

The H\"older moduli $\Theta(\eps)=\eps^\alpha$ are the main case in this
paper. Fernique's theorem gives the integrability of the
random constant in the H\"older modulus of continuity.

Let $(T,d)$ be a compact metric space. For $\alpha>0$, define the H\"older
seminorm on $C(T)$ by
\[
[f]_\alpha
:=
\sup_{\substack{s,t\in T\\s\neq t}}
\frac{|f(s)-f(t)|}{d(s,t)^\alpha}.
\]
It is a measurable, possibly extended-valued seminorm on $C(T)$. For the
isonormal process on a compact set $T\subset H$, define the \emph{critical H\"older
exponent} by
\begin{equation}\label{eq:critical-holder-exponent}
\alpha_H(T)
:=
\sup\Bigl(
\{0\}\cup
\{\beta>0:X\text{ has a }\beta\text{-H\"older modification on }T\}
\Bigr).
\end{equation}

A continuous modification of a centered Gaussian process on a compact metric
space $T$ may be regarded as a $C(T)$-valued random element. For every
$\ell\in C(T)^*$, the random variable $\ell(X)$ is a limit in probability of
finite linear combinations of the coordinates $X_t$ and is therefore
Gaussian. Thus $X$ is a Gaussian random element of $C(T)$.

\begin{theorem}[Fernique's theorem for measurable seminorms
{\cite[Theorem~2.8.5]{Bog98}; see also~\cite{Fer75}}]
Let $Y$ be a centered Gaussian random element of a Banach space $E$, and let
$q\colon E\to[0,\infty]$ be a measurable seminorm. If $q(Y)<\infty$ almost
surely, then there exists $c>0$ such that
\[
\mathbb{E}\exp\bigl(cq(Y)^2\bigr)<\infty.
\]
\end{theorem}

\begin{corollary}\label{cor:fernique-holder}
Let $X$ be a continuous modification of a centered Gaussian process on a
compact metric space $T$. If its sample paths are almost surely
$\alpha$-H\"older, then there exists $c>0$ such that
\[
\mathbb{E}\exp\bigl(c[X]_\alpha^2\bigr)<\infty.
\]
In particular, $\mathbb{E}\,[X]_\alpha<\infty$.
\end{corollary}

\begin{proof}
The preceding discussion shows that $X$ is a Gaussian random element of
$C(T)$, while the H\"older assumption gives $[X]_\alpha<\infty$ almost surely.
The claim follows from Fernique's theorem applied to $q(f)=[f]_\alpha$.
\end{proof}

Thus, if the sample paths are almost surely $\alpha$-H\"older, then
\[
|X_s(\omega)-X_t(\omega)|
\leqslant
C_\alpha(\omega)d(s,t)^\alpha,
\qquad s,t\in T,
\]
almost surely, where the random constant may be chosen as
$C_\alpha=[X]_\alpha$ and satisfies
$\mathbb{E}C_\alpha<\infty$.

\subsection{Local Gaussian width function}
\label{subsec:width-functions}

Let $T\subset H$ be a nonempty compact set. The \emph{local Gaussian width
function} is
\begin{equation}\label{eq:local-gaussian-width}
w_T^{\mathrm{loc}}(\eps)
:=
\sup_{t_0\in T}
\mathbb{E}\sup_{t\in T\cap B(t_0,\eps)}
(X_t-X_{t_0}),
\qquad \eps>0.
\end{equation}

On the one hand, since the local Gaussian width is defined in terms of
expected suprema of Gaussian increments on $T\cap B(t_0,\eps)$, it can be
estimated using Dudley's entropy upper bound~\cite{Dud67} and Sudakov's
minoration (Proposition~\ref{prop:sudakov-minoration}), applied on each such
neighborhood. On the other hand, the following theorem, proved in
Appendix~\ref{app:width-regularity}, relates its small-scale decay to
H\"older regularity of the sample paths.

\begin{theorem}[H\"older regularity and Gaussian width functions]
\label{thm:width-holder-equivalence}
Let $T\subset H$ be a nonempty compact GB-set, and let $0<\beta<1$. The
following conditions are equivalent:
\begin{enumerate}
\item the isonormal process on $T$ has a modification with almost surely
$\beta$-H\"older sample paths;
\item $w_T^{\mathrm{loc}}(\eps)=O(\eps^\beta)$ as $\eps\downarrow0$.
\end{enumerate}
\end{theorem}

\subsection{Intrinsic volumes}
\label{subsec:intrinsic-volumes}

For a nonempty compact convex set $K\subset\mathbb{R}^d$, the intrinsic volumes
are defined by the Steiner formula
\[
\Vol_d(K+\lambda B^d)
=
\sum_{j=0}^d\kappa_{d-j}V_j(K)\lambda^{d-j},
\qquad \lambda\geqslant0,
\]
where $B^d$ is the Euclidean unit ball and
$\kappa_j=\Vol_j(B^j)=\pi^{j/2}/\Gamma(j/2+1)$. If $K$ is full-dimensional,
$V_d(K)$ is its volume and $V_{d-1}(K)$ is half its surface area. In general,
$V_1(K)$ is a constant multiple of the mean width, $V_0(K)=1$, and
$V_j(K)=0$ for $j>\dim K$. This normalization is independent of the ambient
Euclidean space; see, for example,~\cite{Sch14}.

Following Sudakov~\cite{Sud76} and Chevet~\cite{Chev76}, for a compact convex
set $K\subset H$ one therefore defines
\[
V_j(K)
:=
\sup\{V_j(K'):K'\subset K\text{ is compact, convex, and finite-dimensional}\}.
\]
The resulting value may be infinite. Chevet proved that, for every
$j\geqslant1$, one has $V_j(K)<\infty$ if and only if
$V_1(K)<\infty$~\cite{Chev76}. Sudakov's formula
\begin{equation}\label{eq:sudakov-intrinsic-volume}
V_1(K)=\sqrt{2\pi}\,\mathbb{E}\sup_{t\in K}X_t
\end{equation}
shows, in particular, that a compact convex set is GB if and only if its first
intrinsic volume is finite.

The intrinsic volumes satisfy the following dimension-free inequalities, which
follow from the Alexandrov--Fenchel inequality.

\begin{proposition}
\label{prop:intrinsic-volume-inequalities}
Let $K\subset H$ be a nonempty compact convex GB-set. Its intrinsic volumes
satisfy the \emph{ultra-log-concavity} inequality
\begin{equation}\label{eq:intrinsic-ultra-log-concavity}
V_j(K)^2\geqslant\frac{j+1}{j}V_{j-1}(K)V_{j+1}(K),
\qquad j\geqslant1.
\end{equation}
In particular, all intrinsic volumes can be bounded by the first one through
the Chevet--McMullen inequality~\cite{Chev76,McM91}:
\begin{equation}\label{eq:chevet-mcmullen}
V_j(K)\leqslant\frac{V_1(K)^j}{j!},
\qquad j\geqslant1.
\end{equation}
\end{proposition}

Consequently, the generating series of the intrinsic volumes satisfies
\begin{equation}\label{eq:intrinsic-generating-bound}
\sum_{j=0}^{\infty}\lambda^j V_j(K)
\leqslant\exp\bigl(\lambda V_1(K)\bigr)<\infty,
\qquad \lambda\geqslant0.
\end{equation}

Tsirelson's theorem represents every intrinsic volume, up to an explicit
constant, as the expected volume of a Gaussian image. Let $K\subset H$ be a
compact convex GB-set, and let $X^{(1)},\ldots,X^{(j)}$ be independent natural
modifications of the restrictions to $K$ of independent copies of the
isonormal process on $H$. Here naturality means that the sample paths are
almost surely continuous with respect to some separable metric on $K$; this
stronger choice of modification is essential for the spectrum construction.
Define the random set
\[
\operatorname{Spec}_j K
:=
\overline{\left\{
\bigl(X_t^{(1)},\ldots,X_t^{(j)}\bigr):t\in K
\right\}}
\subset\mathbb{R}^j.
\]
For these natural modifications, Tsirelson's theorem implies that the
image is almost surely bounded and convex~\cite[Theorem~6]{Tsi85}. Thus
$\operatorname{Spec}_j K$ is compact and convex, and taking the closure does not
change its $j$-dimensional volume. This set is called the $j$-dimensional
\emph{Gaussian spectrum} of $K$.

\begin{theorem}[Tsirelson~{\cite[Theorem~6]{Tsi85}}]
\label{thm:tsirelson-projection}
For every compact convex GB-set $K\subset H$ and every $j\geqslant1$,
\[
V_j(K)
=
\frac{(2\pi)^{j/2}}{j!\kappa_j}
\,\mathbb{E}\Vol_j\bigl(\operatorname{Spec}_j K\bigr),
\qquad
\kappa_j=\Vol_j(B^j).
\]
\end{theorem}

For $j=1$, the theorem reduces to~\eqref{eq:sudakov-intrinsic-volume}.

Tsirelson also gave a Gaussian representation of the generating function
of the intrinsic volumes. This is the infinite-dimensional Steiner
formula recalled in~\eqref{eq:tsirelson-steiner}.

\begin{theorem}[Tsirelson~{\cite[Theorem~5]{Tsi85}}]
\label{thm:tsirelson-steiner}
Let $K\subset H$ be a nonempty compact convex GB-set, and let
$(X_t)_{t\in H}$ be the isonormal process on $H$. Then, for every
$\lambda\geqslant0$,
\[
\mathbb{E}\exp\Bigl(\sup_{t\in K}
\bigl[\sqrt{2\pi}\,\lambda X_t-\pi\lambda^2\|t\|^2\bigr]\Bigr)
=\sum_{j=0}^{\infty}\lambda^j V_j(K)<\infty.
\]
The expectation on the left is understood as the supremum of the
corresponding expectations over nonempty finite subsets of $K$.
\end{theorem}

\subsection{The Gaussian Wills functional}

The Gaussian expression on the left-hand side of
Theorem~\ref{thm:tsirelson-steiner} makes sense for arbitrary
nonempty index sets. Accordingly, for a nonempty set $T\subset H$ we define
the \emph{Gaussian Wills functional}
\begin{equation}\label{eq:wills-gaussian}
W(T):=\mathbb{E}\exp\Bigl(\sup_{t\in T}\bigl[\sqrt{2\pi}\,X_t-\pi\|t\|^2\bigr]\Bigr)
:=\sup_{T_0\subset T\text{ finite}}
\mathbb{E}\exp\Bigl(\max_{t\in T_0}\bigl[\sqrt{2\pi}\,X_t-\pi\|t\|^2\bigr]\Bigr),
\end{equation}
where, just as for the expected supremum, the expectation is defined as the
supremum over finite subsets of $T$. Thus $W(T)\in[1,\infty]$.
For every compact convex GB-set $K\subset H$, Tsirelson's
formula~\eqref{eq:tsirelson-steiner} gives
$W(\lambda K)=\sum_{j=0}^{\infty}\lambda^jV_j(K)$ for every
$\lambda\geqslant0$. This definition agrees with the intrinsic-volume
definition on convex compact sets and extends it to nonconvex index sets.

We next derive the analogue of the integral
representation~\eqref{eq:integral-repr-wills} for
arbitrary finite-dimensional sets. Mourtada~\cite{Mou25}
was the first to consider the Wills functional in this integral form for
sets that are not necessarily convex. This representation makes the metric
nature of the functional clear. Indeed, let $T\subset H$ be a
finite-dimensional set, and identify a
finite-dimensional linear subspace containing $T$ with $\mathbb{R}^d$. Then
$(X_t)$ becomes
$(\la Z_d,t\ra)$, where $Z_d$ is a standard Gaussian vector in
$\mathbb{R}^d$.
Using
$\sqrt{2\pi}\la z,t\ra-\pi\|t\|^2
=\tfrac12\bigl(\|z\|^2-\|z-\sqrt{2\pi}\,t\|^2\bigr)$, we get
\[
W(T)
=\mathbb{E}\exp\Bigl(\sup_{t\in T}
\bigl[\sqrt{2\pi}\la Z_d,t\ra-\pi\|t\|^2\bigr]\Bigr)
=\mathbb{E}\exp\Bigl(\frac{\|Z_d\|^2-\dist^2(Z_d,\sqrt{2\pi}\,T)}{2}\Bigr).
\]
When the expectation is written as an integral against the standard Gaussian
density, the factor $e^{\|z\|^2/2}$ cancels, and the substitution
$z=\sqrt{2\pi}\,x$ gives
\begin{equation}\label{eq:wills-integral}
W(T)
=\frac{1}{(2\pi)^{d/2}}\int_{\mathbb{R}^d}
\exp\Bigl(-\frac{\dist^2(z,\sqrt{2\pi}\,T)}{2}\Bigr)dz
=\int_{\mathbb{R}^d}\exp\left(-\pi\dist^2(x,T)\right)dx.
\end{equation}

The right-hand side of~\eqref{eq:wills-integral} depends only on the metric
space $(T,\|\cdot\|)$. Moreover, it does not depend on the dimension of the
ambient space. Indeed, for $D>d$ embed $\mathbb{R}^d$ into
$\mathbb{R}^D=\mathbb{R}^d\times\mathbb{R}^{D-d}$. Then
$\dist^2((x,y),T)=\dist^2(x,T)+\|y\|^2$ for $x\in\mathbb{R}^d$ and
$y\in\mathbb{R}^{D-d}$, and since $e^{-\pi\|y\|^2}$ is a probability density
on $\mathbb{R}^{D-d}$, Fubini's theorem gives
\[
\int_{\mathbb{R}^D}e^{-\pi\dist^2((x,y),T)}\,dx\,dy
=\int_{\mathbb{R}^d}e^{-\pi\dist^2(x,T)}\,dx
\int_{\mathbb{R}^{D-d}}e^{-\pi\|y\|^2}\,dy
=\int_{\mathbb{R}^d}e^{-\pi\dist^2(x,T)}\,dx.
\]

\begin{definition}\label{def:wills-growth}
Let $T\subset H$ be a nonempty compact set. If $T$ is not a GB-set, we set
$\rho(T):=+\infty$. Otherwise, its \emph{Wills growth exponent} is
\[\rho(T):=\limsup_{r\to\infty}\frac{\log\log W(rT)}{\log r}.\]
As usual, a constant function $r\mapsto W(rT)$ is assigned exponent zero.
\end{definition}

We next collect the properties of $W$ used throughout.

\begin{proposition}[Elementary properties of the Wills functional]
  \label{prop:wills-properties}
  Let $A, B, A_1, \ldots, A_n$ be nonempty subsets of $H$.
  \begin{enumerate}
    \item \emph{(Monotonicity.)} If $A \subset B$, then $W(A) \leqslant W(B)$.
    \item \emph{(Isometry invariance.)} If $\varphi \colon A \to B$ is a
      surjective isometry, then $W(A) = W(B)$.
    \item \emph{(Subadditivity.)} If $A \subset \bigcup_{i=1}^{n} A_i$, then
      $W(A) \leqslant \sum_{i=1}^{n} W(A_i)$.
    \item \emph{(Exponential bound.)} If $\mathbb{E}\sup_{t \in A}X_t < \infty$,
      then for every $r > 0$,
      \[\log W(rA) \leqslant \sqrt{2\pi}\,r\,\mathbb{E}\sup_{t \in A}X_t.\]
  \end{enumerate}
\end{proposition}

\begin{proof}
  Property~(1) is immediate from~\eqref{eq:wills-gaussian}, since every finite
  subset of $A$ is a finite subset of $B$.

  Property~(2) follows from the integral representation~\eqref{eq:wills-integral}
  applied to finite-dimensional subsets.

  For~(3), put $Y_t:=\sqrt{2\pi}X_t-\pi\|t\|^2$. It suffices to prove
  the inequality when all the sets are finite. Since
  $A\subset\bigcup_{i=1}^n A_i$ and
  $\exp(\max_i y_i)\leqslant\sum_i\exp(y_i)$ for real numbers $y_i$,
  \[
  W(A)=\mathbb{E}\exp\Bigl(\max_{t\in A}Y_t\Bigr)
  \leqslant\sum_{i=1}^n\mathbb{E}\exp\Bigl(\max_{t\in A_i}Y_t\Bigr)
  =\sum_{i=1}^n W(A_i).
  \]
  For arbitrary sets, taking suprema over finite subsets
  proves~(3).

  For~(4), fix a finite
  set $A_0\subset A$ and let $K_0=\conv A_0$. Then $K_0$ is a finite-dimensional
  compact convex set, the supremum of the linear function $t\mapsto X_t$ over
  $K_0$ equals its maximum over $A_0$, and Sudakov's formula gives
  $V_1(K_0)=\sqrt{2\pi}\,\mathbb{E}\max_{t\in A_0}X_t$. By monotonicity, the
  Chevet--McMullen bound~\cite{Chev76,McM91}, and the homogeneity
  $V_1(rK_0)=rV_1(K_0)$,
  \[W(rA_0)\leqslant W(rK_0)\leqslant \exp\left(rV_1(K_0)\right)
  \leqslant \exp\left(\sqrt{2\pi}\,r\,\mathbb{E}\sup_{t\in A}X_t\right).\]
  Taking the supremum over finite $A_0\subset A$ proves the claim.
\end{proof}

\begin{lemma}\label{lem:wills-finite-gb}
Let $T\subset H$ be a nonempty compact set. Then $T$ is Gaussian bounded if and
only if $W(rT)<\infty$ for some $r>0$; in this case $W(rT)<\infty$ for all
$r>0$.
\end{lemma}

\begin{proof}
By Sudakov's criterion recalled in the Introduction, $T$ is Gaussian bounded if
and only if $\mathbb{E}\sup_{t\in T}X_t<\infty$. If $T$ is Gaussian bounded, the
exponential bound of Proposition~\ref{prop:wills-properties} gives
$W(rT)<\infty$ for all $r>0$. Conversely, assume $W(r_0T)<\infty$ for some
$r_0>0$, and put $R:=\sup_{t\in T}\|t\|<\infty$. For every finite
$T_0\subset T$, Jensen's inequality gives
\[\log W(r_0T)
\geqslant
\log\mathbb{E}\exp\left(
\max_{t\in T_0}
\left[\sqrt{2\pi}r_0X_t-\pi r_0^2\|t\|^2\right]\right)
\geqslant
\sqrt{2\pi}r_0\,\mathbb{E}\max_{t\in T_0}X_t
-\pi r_0^2R^2 .\]
Taking the supremum over all finite $T_0\subset T$ shows
$\mathbb{E}\sup_{t\in T}X_t<\infty$, so the same criterion shows that $T$ is
Gaussian bounded.
\end{proof}

Finally, we recall Mourtada's comparison principle, the Wills-functional analogue
of the Sudakov--Fernique comparison principle: a contraction of index sets
does not increase the Wills functional, just as it does not increase the
expected supremum.

\begin{proposition}[Mourtada's comparison principle~{\cite[Theorem~4.1]{Mou25}}]
  \label{prop:comparison}
  Let $A$ and $B$ be nonempty subsets of separable Hilbert spaces. If there
  exists a contraction $\varphi \colon A \to B$, i.e., a map satisfying
  $\left\|\varphi(x) - \varphi(y)\right\| \leqslant \left\|x - y\right\|$
  for all $x,y \in A$, with $\varphi(A) = B$, then $W(A) \geqslant W(B)$.
\end{proposition}

%% file: main_results.tex
\section{Main results}\label{sec:main}

\subsection{From Wills growth to metric entropy}

The first step is the Wills-functional analogue of the elementary lemma used in
the proof of Sudakov's minoration. An $\eps$-separated finite set can be
contracted onto the vertices of a regular simplex with mutual distance $\eps$;
Mourtada's comparison principle then reduces the estimate to an
explicit computation for this discrete simplex.

\begin{lemma}\label{lem:entropy-lower-bound}
  There exist universal constants $c_1, c_2 > 0$ such that, for every
  $\eps$-separated set $\left\{\theta_1, \ldots, \theta_N\right\} \subset H$,
  \[\log W\left(\left\{\theta_1, \ldots, \theta_N\right\}\right) \geqslant c_1 \eps \sqrt{\log N} - c_2 \eps^2.\]
\end{lemma}

\begin{proof}
  Let $e_1, \ldots, e_N$ be the standard
  orthonormal basis of $\mathbb{R}^N$, and consider the vertex set of a regular
  simplex scaled to edge length $\eps$:
  \[S = \left\{v_i = \frac{\eps}{\sqrt{2}} e_i : 1 \leqslant i \leqslant N\right\},
  \qquad \left\|v_i - v_j\right\| = \eps \quad (i \neq j).\]
  Since $\left\|\theta_i - \theta_j\right\| \geqslant \eps$ for all $i \neq j$,
  the map $\varphi\left(\theta_i\right) = v_i$ is a contraction from
  $\left\{\theta_1, \ldots, \theta_N\right\}$ onto $S$, and the comparison
  principle (Proposition~\ref{prop:comparison}) gives
  \[W\left(\left\{\theta_1, \ldots, \theta_N\right\}\right) \geqslant W(S).\]
  Let $Z=(Z_1,\ldots,Z_N)\sim\mathcal{N}\left(0,I_N\right)$. Since
  $\left\|v_i\right\|^2 = \frac{\eps^2}{2}$ for all $i$, the Gaussian
  representation gives
  \begin{align*}
  W(S)
  &= \mathbb{E} \exp\left(\max_{1 \leqslant i \leqslant N}
     \left[\sqrt{2\pi}\left\langle v_i, Z\right\rangle - \pi \left\|v_i\right\|^2\right]\right)\\
  &= \exp\left(-\frac{\pi \eps^2}{2}\right) \mathbb{E}
     \exp\left(\sqrt{\pi} \eps \max_{1 \leqslant i \leqslant N} Z_i\right).
  \end{align*}
  Jensen's inequality $\mathbb{E} e^X \geqslant e^{\mathbb{E} X}$ and the
  standard lower bound
  \[\mathbb{E} \max_{1 \leqslant i \leqslant N} Z_i \geqslant c' \sqrt{\log N},\]
  valid for $N\geqslant2$ with a universal $c'>0$, yield
  \[W(S) \geqslant \exp\left(c_1 \eps \sqrt{\log N} - c_2 \eps^2\right)\]
  with $c_1 = \sqrt{\pi} c'$ and $c_2 = \frac{\pi}{2}$.
\end{proof}

\begin{theorem}[Wills--Sudakov minoration]\label{thm:w_lower_bound}
  Let $T$ be a nonempty compact set in $H$, and let $c_1, c_2 > 0$ be the
  universal constants of Lemma~\ref{lem:entropy-lower-bound}. Then for all
  $r,\eps > 0$,
  \[\log W(rT) \geqslant c_1 r \eps \sqrt{\log M(T,\eps)} - c_2 r^2 \eps^2,\]
  and consequently, by \eqref{eq:packing-covering}, the same lower bound holds
  with $M(T,\eps)$ replaced by $N(T,\eps)$.
\end{theorem}

\begin{proof}
  Put $M = M(T,\eps)$ and choose an $\eps$-separated set
  $\left\{\theta_1, \ldots, \theta_M\right\} \subset T$ of cardinality $M$. Then
  $\left\{r\theta_1, \ldots, r\theta_M\right\}\subset rT$ is $r\eps$-separated
  (cf.~the scaling relation \eqref{eq:scaling}), and by monotonicity of $W$ and Lemma~\ref{lem:entropy-lower-bound},
  \[\log W(rT) \geqslant \log W\left(\left\{r\theta_1, \ldots, r\theta_M\right\}\right) \geqslant c_1 r \eps \sqrt{\log M} - c_2 r^2 \eps^2.\]
  The version with $N(T,\eps)$ follows since
  $x \mapsto c_1 r\eps\sqrt{x} - c_2 r^2\eps^2$ is increasing.
\end{proof}

Optimizing over the dilation parameter $r$ yields an upper bound on the metric
entropy.

\begin{theorem}[Theorem~B]\label{thm:entropy-from-wills}
Let $T\subset H$ be a nonempty compact set. Then there exists a universal
constant $C>0$ such that, for every $\eps>0$,
\[\sqrt{\log N(T,\eps)}
\leqslant
C\inf_{r>0}
\left(
\frac{\log W(rT)}{r\eps}+r\eps
\right).\]
\end{theorem}

\begin{proof}
By Theorem~\ref{thm:w_lower_bound}, for every $r>0$,
\[c_1r\eps\sqrt{\log N(T,\eps)}
\leqslant
\log W(rT)+c_2r^2\eps^2.\]
Dividing by $c_1 r\eps$ and taking the infimum over $r>0$ proves the claim with
$C=\max\left(c_1^{-1},c_2c_1^{-1}\right)$.
\end{proof}

Now, assuming power-law bounds on $\log W(rT)$, one can obtain power-law bounds on $\log N(T,\varepsilon)$, and hence on the Dudley integral. Below we state the corollaries of Theorem B, with proofs given in Appendix~\ref{app:proofs}.

\begin{corollary}\label{cor:entropy-upper-bound}
Let $T\subset H$ be a nonempty compact set with $\rho(T)<1$. Then, for every
\[
0<\alpha<\frac{2-2\rho(T)}{2-\rho(T)},
\]
there exists a constant $C_{T,\alpha}>0$ such that
\[\log N(T,\eps)
\leqslant
C_{T,\alpha}
\eps^{-(2-2\alpha)}\]
for all $0<\eps\leqslant \diam T$.
\end{corollary}

\begin{corollary}\label{cor:holder}
Let $T\subset H$ be a nonempty compact set with $\rho(T)<1$, and let
$(X_t)_{t\in T}$ be the isonormal process indexed by $T$. Then, for every
\[
0<\alpha<\frac{2-2\rho(T)}{2-\rho(T)},
\]
there exist a modification $(\widetilde X_t)_{t\in T}$ and a deterministic
constant $C_{T,\alpha}<\infty$ such that
\[\limsup_{\eps\to0}
\frac{\theta(\widetilde X(\omega),\eps)}
{\eps^\alpha}
\leqslant C_{T,\alpha}\]
almost surely. In particular, this modification has almost surely
$\alpha$-H\"older sample paths on $T$.
\end{corollary}

% % \begin{remark}\label{rem:holder-uniform}
% %   The modulus in Corollary~\ref{cor:holder} is uniform, not merely local: the
% %   supremum in $\theta$ runs over all pairs $s,t\in T$ with
% %   $\|s-t\|\leqslant\eps$. The critical exponent $\frac{2-2\rho(T)}{2-\rho(T)}$
% %   is approached from below; the corollary does not claim the endpoint exponent,
% %   and Example~\ref{ex:wiener-spiral} shows that in general the endpoint fails.
% % \end{remark}

% \begin{remark}\label{rem:boundary-case}
% The statement gives nontrivial H\"older exponents only when $\rho(T)<1$; at
% $\rho(T)=1$ the exponent degenerates to zero. This does not mean that
% continuity is impossible in the boundary case: for convex $K$ with
% $\rho(K)=1$, Gaussian continuity is controlled by the type of the Steiner
% entire function \cite{DGZ25}, and the boundary case naturally leads to
% non-H\"older moduli with logarithmic corrections. The quantitative estimate of
% Theorem~\ref{thm:entropy-from-wills} still applies whenever more precise
% information on the growth of $W(rT)$ is available.
% \end{remark}

\subsection{From local Gaussian widths to Wills growth}

We estimate the Wills functional in terms of local Gaussian widths (see \eqref{eq:local-gaussian-width}). The key
idea is to apply subadditivity to a covering by $\eps$-balls; this argument goes
back to Vitale~\cite{Vit01}.

\begin{proposition}\label{prop:wills-covering}
For every nonempty compact set $T\subset H$ and all $r,\eps>0$,
\[
\log W(rT)
\leqslant
\log N(T,\eps)+\sqrt{2\pi}\,r w_T^{\mathrm{loc}}(\eps).
\]
\end{proposition}

\begin{proof}
If $w_T^{\mathrm{loc}}(\eps)=\infty$, there is nothing to prove. Otherwise,
let $t_1,\ldots,t_{N(T,\eps)}\in T$ be the centers of a minimal
$\eps$-cover and put $T_i=T\cap B(t_i,\eps)$. Subadditivity, isometry
invariance, and the exponential bound in Proposition~\ref{prop:wills-properties}
give
\[
W(rT)
\leqslant
\sum_i W(rT_i)
=\sum_i W\bigl(r(T_i-t_i)\bigr)
\leqslant
N(T,\eps)\exp\left(\sqrt{2\pi}\,r
w_T^{\mathrm{loc}}(\eps)\right).
\]
Here we used
$\mathbb{E}\sup_{t\in T_i}(X_t-X_{t_i})
\leqslant w_T^{\mathrm{loc}}(\eps)$. Taking the logarithm finishes the proof.
\end{proof}

We next estimate the metric entropy. Sudakov's minoration inside each ball
bounds $\log N(T,\delta)-\log N(T,2\delta)$, and summing these differences gives the
following result.

\begin{lemma}\label{lem:entropy-from-local-width}
There exists a universal constant $C>0$ such that, for every nonempty compact
set $T\subset H$ and every $\eps>0$,
\[
\log N(T,\eps)
\leqslant
C\int_{\eps}^{\infty}
\frac{w_T^{\mathrm{loc}}(t)^2}{t^3}\,dt.
\]
\end{lemma}

\begin{proof}
Let $D=\diam T$. The claim is immediate for $\eps\geqslant D$, so assume
$0<\eps<D$. For $\delta>0$, if
$t_1,\ldots,t_{N(T,2\delta)}$ are the centers of a minimal $2\delta$-cover and
$T_i=T\cap B(t_i,2\delta)$, then
$N(T,\delta)\leqslant N(T,2\delta)\max_i N(T_i,\delta)$. Sudakov's minoration inside each
$T_i$ gives
\[
\log N(T,\delta)-\log N(T,2\delta)
\leqslant
C\frac{w_T^{\mathrm{loc}}(2\delta)^2}{\delta^2}.
\]
Indeed,
\[
\delta\sqrt{\log N(T_i,\delta)}
\leqslant
C\mathbb{E}\sup_{t\in T_i}(X_t-X_{t_i})
\leqslant Cw_T^{\mathrm{loc}}(2\delta).
\]
Summing over $\delta=2^k\eps$ until $2^k\eps\geqslant D$ gives
\[
\log N(T,\eps)
\leqslant
C\sum_{k=0}^{\infty}
\frac{w_T^{\mathrm{loc}}(2^{k+1}\eps)^2}{2^{2k}\eps^2}.
\]
Since $T\cap B(t_0,\delta)$ increases with $\delta$, the function
$w_T^{\mathrm{loc}}$ is nondecreasing. Put $s_k=2^{k+1}\eps$. Then
\[
\frac{w_T^{\mathrm{loc}}(2^{k+1}\eps)^2}{2^{2k}\eps^2}
\leqslant
\frac{32}{3}
\int_{s_k}^{2s_k}\frac{w_T^{\mathrm{loc}}(t)^2}{t^3}\,dt.
\]
The intervals $[s_k,2s_k]$ are disjoint, so summation proves the claim.
\end{proof}

We are now ready to bound the Wills functional in terms of the local Gaussian width.
\begin{theorem}[Theorem~C]
\label{thm:wills-from-local-width}
There exists a universal constant $C>0$ such that, for every nonempty compact
$T\subset H$ and every $r>0$,
\[
\log W(rT)
\leqslant
\inf_{\eps>0}\left\{
C\int_{\eps}^{\infty}
\frac{w_T^{\mathrm{loc}}(t)^2}{t^3}\,dt
+\sqrt{2\pi}\,r w_T^{\mathrm{loc}}(\eps)
\right\}.
\]
\end{theorem}

\begin{proof}
Proposition~\ref{prop:wills-covering} and
Lemma~\ref{lem:entropy-from-local-width} give, for every $\eps>0$,
\[
\log W(rT)
\leqslant
C\int_{\eps}^{\infty}
\frac{w_T^{\mathrm{loc}}(t)^2}{t^3}\,dt
+\sqrt{2\pi}\,r w_T^{\mathrm{loc}}(\eps),
\]
and taking the infimum proves the claim.
\end{proof}

Suppose that $w_T^{\mathrm{loc}}(\eps)\leqslant C\eps^\beta$ for small $\eps$.
Since $w_T^{\mathrm{loc}}$ is bounded on a compact GB-set, the same inequality
holds for all $\eps>0$ with a possibly larger constant $C$. Substitution in Theorem~C gives
the following result; its proof is in Appendix~\ref{app:proofs}.

\begin{corollary}\label{cor:order-from-local-width}
Let $T\subset H$ be a nonempty compact GB-set. Suppose that, for some
$\beta\in(0,1)$, $C>0$, and $\eps_0>0$,
\[
w_T^{\mathrm{loc}}(\eps)\leqslant C\eps^\beta
\]
for every $0<\eps\leqslant\eps_0$. Then there exist constants $C'>0$ and
$r_0\geqslant1$ such that, for every $r\geqslant r_0$,
\[
\log W(rT)\leqslant
C' r^{\frac{2-2\beta}{2-\beta}}.
\]
In particular,
\[
\rho(T)\leqslant\frac{2-2\beta}{2-\beta}.
\]
\end{corollary}

\subsection{Relations between the exponents}

Combining Theorems~B and~C with Dudley's theorem yields the explicit relations
between the entropy exponent, the critical H\"older exponent, and the Wills
growth exponent.

\begin{theorem}[Theorem~A]
\label{thm:asymptotic-exponents}
Let $T\subset H$ be a compact GB-set. Then
\begin{equation}\label{eq:entropy-wills-identity}
\eta(T)
=
\frac{2\rho(T)}{2-\rho(T)},
\qquad
\rho(T)
=
\frac{2\eta(T)}{2+\eta(T)},
\end{equation}
and
\begin{equation}\label{eq:order-holder-identity}
\alpha_H(T)
=
1-\frac{\eta(T)}{2}
=
\frac{2-2\rho(T)}{2-\rho(T)}.
\end{equation}
\end{theorem}

\begin{proof}
Write $\alpha=\alpha_H(T)$, $\eta=\eta(T)$, and $\rho=\rho(T)$. We prove
\[
\alpha
\geqslant 1-\frac{\eta}{2}
\geqslant \frac{2-2\rho}{2-\rho}
\geqslant \alpha.
\]

If $\eta<2$, then for every $\eta'\in(\eta,2)$ the definition of $\eta$ gives
$\log N(T,\eps)=O(\eps^{-\eta'})$. Hence
$\mathcal D(T,u)=O(u^{1-\eta'/2})$, and Theorem~\ref{thm:dudley}, followed by
compactness, gives a $(1-\eta'/2)$-H\"older modification. Letting
$\eta'\downarrow\eta$ proves the first inequality; for $\eta=2$ it is automatic.

If $\rho<1$, then for every $\rho'\in(\rho,1)$ we have
$\log W(rT)=O(r^{\rho'})$. Corrollary ~\ref{cor:entropy-upper-bound} gives
\[
\log N(T,\eps)=O\bigl(\eps^{-2\rho'/(2-\rho')}\bigr).
\]
Letting $\rho'\downarrow\rho$ yields $\eta\leqslant2\rho/(2-\rho)$, which is
equivalent to the second inequality. For $\rho=1$, it follows from
$\eta\leqslant2$ by Proposition~\ref{prop:sudakov-minoration}.

Finally, for every $0<\alpha'<\alpha$,
Theorem~\ref{thm:width-holder-equivalence} gives
$w_T^{\mathrm{loc}}(\eps)=O(\eps^{\alpha'})$.
By Corollary~\ref{cor:order-from-local-width},
\[
\rho\leqslant\frac{2-2\alpha'}{2-\alpha'}.
\]
Letting $\alpha'\uparrow\alpha$ and rearranging proves the third inequality.
For $\alpha=0$, it follows from $\rho\leqslant1$.

Thus all three inequalities are equalities, and the asserted identities follow.
\end{proof}

\subsection{Projection volumes and the Milman--Pisier theorem}

We conclude with an application of Theorem~A to maximal projection volumes, which, for a compact convex set $K\subset H$, are defined as
\[
\mathcal V_n(K)
:=
\sup_{\operatorname{rank}P=n}\Vol_n(PK),
\]
where the supremum runs over all rank-$n$ orthogonal projections $P$ on $H$,
and $\Vol_n$ denotes Lebesgue measure on the $n$-dimensional
Euclidean space $PH$.

Following Dudley, define the \emph{volume
exponent}
\begin{equation}\label{eq:volume-exponent}
\operatorname{EV}(K)
:=
\limsup_{n\to\infty}
\frac{\log\mathcal V_n(K)}{n\log n}.
\end{equation}

By Tsirelson's
formula, $V_n(K)$ is proportional to the mean volume of a Gaussian image; see
Theorem~\ref{thm:tsirelson-projection}. In contrast,
$\mathcal V_n(K)$ is the maximal volume of a deterministic orthogonal
projection. The theorem of Milman and Pisier relates the latter volumes to
metric entropy.

\begin{theorem}[Milman--Pisier {\cite[Theorem~1]{MP87}}]
\label{thm:milman-pisier}
Let $K=-K$ be a compact convex subset of $H$. Suppose that there
exist $C<\infty$ and $\delta>0$ such that
\begin{equation}\label{eq:milman-pisier-assumption}
\mathcal V_n(K)^{1/n}\leqslant Cn^{-1-\delta}
\end{equation}
for every $n\geqslant1$. Then $K$ is a GC-set. Moreover,
\begin{equation}\label{eq:milman-pisier-identity}
\operatorname{EV}(K)+\frac12=-\frac1{\eta(K)}.
\end{equation}
\end{theorem}

We now assume that $0<\rho(K)<1$. Theorem~A gives the entropy estimate used
below to verify the hypothesis~\eqref{eq:milman-pisier-assumption} of
Theorem~\ref{thm:milman-pisier} for some $\delta>0$.

\begin{corollary}[Maximal projection volumes and Wills growth]
\label{cor:milman-pisier-holder}
Let $K=-K$ be a compact convex GB-set such that
$0<\rho(K)<1$. Then
\begin{equation}\label{eq:volume-wills-identity}
\operatorname{EV}(K)=-\frac1{\rho(K)}.
\end{equation}
\end{corollary}

\begin{proof}
By Theorem~\ref{thm:asymptotic-exponents},
\[
\eta(K)=\frac{2\rho(K)}{2-\rho(K)}<2.
\]
Fix $q$ with $\eta(K)<q<2$. Then
$\log N(K,\eps)\leqslant\eps^{-q}$ for all sufficiently small $\eps$.
Writing $\kappa_n$ for the volume of the Euclidean unit ball in
$\mathbb R^n$, projecting an $\eps$-cover and comparing volumes gives
\[
\mathcal V_n(K)^{1/n}
\leqslant
\kappa_n^{1/n}\eps N(K,\eps)^{1/n}.
\]
Set $\eps=n^{-1/q}$. For all sufficiently large $n$,
$N(K,\eps)^{1/n}\leqslant e$, and $\kappa_n^{1/n}=O(n^{-1/2})$; hence
\[
\mathcal V_n(K)^{1/n}
\leqslant
C_q n^{-1/2-1/q}
=C_q n^{-1-\delta_q},
\qquad
\delta_q:=\frac1q-\frac12>0.
\]
After increasing $C_q$ for the finitely many remaining $n$,
Theorem~\ref{thm:milman-pisier} applies. Its identity and Theorem~A give
\[
\operatorname{EV}(K)
=-\frac1{\eta(K)}-\frac12
=-\frac{2-\rho(K)}{2\rho(K)}-\frac12
=-\frac1{\rho(K)}.
\]
\end{proof}

\begin{remark}[Intrinsic volumes versus maximal projection volumes]
\label{rem:tsirelson-milman-pisier}
The coefficient formula \eqref{eq:order-from-coefficients}, applied to the
Steiner entire function, gives
\[
\limsup_{n\to\infty}
\frac{\log V_n(K)}{n\log n}
=
-\frac1{\rho(K)}.
\]
Consequently, under the assumptions of
Corollary~\ref{cor:milman-pisier-holder},
\[
\limsup_{n\to\infty}
\frac{\log V_n(K)}{n\log n}
=
\operatorname{EV}(K) = \limsup_{n\to\infty}
\frac{\log\mathcal V_n(K)}{n\log n}.
\]
\end{remark}

%% file: appendix.tex
\appendix
\section{Proofs of the Corollaries}\label{app:proofs}

\begin{proof}[Proof of Corollary~\ref{cor:entropy-upper-bound}]
Fix
\[
0<\alpha<\frac{2-2\rho(T)}{2-\rho(T)}.
\]
Choose $\beta\in(\rho(T),1)$ sufficiently close to $\rho(T)$ so that
\[
\alpha<\frac{2-2\beta}{2-\beta}.
\]
The
definition of $\rho(T)$ then implies that there exist
$A>0$ and $r_0\geqslant1$ such that
\[\log W(rT)\leqslant Ar^\beta,\qquad r\geqslant r_0.\]
For $\eps>0$ small enough, the choice
\[r_\eps=A^{1/(2-\beta)}\eps^{-2/(2-\beta)}\]
satisfies $r_\eps\geqslant r_0$, and Theorem~\ref{thm:entropy-from-wills} gives
\[\sqrt{\log N(T,\eps)}
\leqslant
C\left(
\frac{A r_\eps^\beta}{r_\eps\eps}
+r_\eps\eps
\right)
\leqslant
C_\beta'\eps^{-\beta/(2-\beta)}.\]
Squaring gives
\[
\log N(T,\eps)\leqslant C_\beta''\eps^{-\frac{2\beta}{2-\beta}}
\]
for all sufficiently small $\eps>0$. The choice of $\beta$ gives
\[
\frac{2\beta}{2-\beta}<2-2\alpha,
\]
so the claimed estimate follows for all sufficiently small $\eps$. Since the covering numbers are finite and nonincreasing in $\eps$, we can increase the constant to cover the remaining values of $\eps$ as well. Thus, taking $C_{T,\alpha}$ possibly larger than $C_\beta''$, we obtain
\[
\log N(T,\eps)
\leqslant
C_{T,\alpha}\eps^{-(2-2\alpha)}
\]
for all $0<\eps\leqslant\diam T$.
\end{proof}

\begin{proof}[Proof of Corollary~\ref{cor:holder}]
If $T$ is finite, take $\widetilde X=X$. Then
$\theta(\widetilde X(\omega),u)=0$ for all sufficiently small $u>0$, so
the conclusion is immediate. Thus, assume that $T$ is infinite. Fix
\[
0<\alpha<\frac{2-2\rho(T)}{2-\rho(T)}.
\]
By Corollary~\ref{cor:entropy-upper-bound},
$\log N(T,\eps)\leqslant A_{T,\alpha}\eps^{-(2-2\alpha)}$ for all sufficiently small
$\eps>0$, whence
\[\mathcal{D}(T,u)
=
\int_0^u\sqrt{\log N(T,\eps)}\,d\eps
\leqslant
\sqrt{A_{T,\alpha}}
\int_0^u\eps^{-(1-\alpha)}\,d\eps
=
\frac{\sqrt{A_{T,\alpha}}}{\alpha}
u^\alpha.\]
In particular, the Dudley integral is finite, so by Theorem~\ref{thm:dudley}
the process has an almost surely uniformly continuous modification
$\widetilde X$ satisfying
\[\limsup_{u\to0}
\frac{\theta(\widetilde X(\omega),u)}{u^\alpha}
\leqslant
\frac{36\sqrt{A_{T,\alpha}}}{\alpha}
=:C_{T,\alpha}\]
almost surely.

The estimate above shows that for almost every $\omega$ there exists
$\delta(\omega)>0$ such that
$|\widetilde X_s(\omega)-\widetilde X_t(\omega)|\leqslant(C_{T,\alpha}+1)\|s-t\|^\alpha$
whenever $\|s-t\|\leqslant\delta(\omega)$. Since $T$ is
compact, each continuous sample path is bounded, so allowing a larger random
constant extends this bound to all pairs of points in $T$. Thus the
modification is $\alpha$-H\"older on all of $T$.
\end{proof}

\begin{proof}[Proof of Corollary~\ref{cor:order-from-local-width}]
Choose $\delta>0$ so that the assumed bound holds for
$0<t\leqslant\delta$, and put
$M:=\mathbb{E}\sup_{t\in T}X_t$, which is finite.
The process is centered, so $w_T^{\mathrm{loc}}(t)\leqslant M$ for every $t>0$.
Since $2\beta-2<0$, for
$0<\eps\leqslant\delta$ we have
\[
\int_{\eps}^{\infty}
\frac{w_T^{\mathrm{loc}}(t)^2}{t^3}\,dt
\leqslant
C^2\int_{\eps}^{\delta}t^{2\beta-3}\,dt
+M^2\int_{\delta}^{\infty}t^{-3}\,dt
\leqslant A\eps^{2\beta-2}.
\]
Theorem~\ref{thm:wills-from-local-width} therefore gives
\[
\log W(rT)\leqslant A\eps^{2\beta-2}+Br\eps^\beta.
\]
For all sufficiently large $r$, the choice
$\eps=r^{-1/(2-\beta)}$ lies below $\delta$ and balances the two terms, giving
\[
\log W(rT)\leqslant
(A+B)r^{\frac{2-2\beta}{2-\beta}}.
\]
Taking logarithms once more, dividing by $\log r$, and letting $r\to\infty$
yields the asserted bound on $\rho(T)$.
\end{proof}

\section{Critical H\"older regularity and Gaussian width functions}
\label{app:width-regularity}

We prove Theorem~\ref{thm:width-holder-equivalence}, which relates the
local Gaussian width to the H\"older exponent.

\begin{proof}[Proof of Theorem~\ref{thm:width-holder-equivalence}]
Assume first that the process has a $\beta$-H\"older modification and fix this
modification. Corollary~\ref{cor:fernique-holder} gives
$\mathbb{E}[X]_\beta<\infty$. Consequently,
\[
w_T^{\mathrm{loc}}(\eps)
\leqslant
\mathbb{E}[X]_\beta\,\eps^\beta,
\]
which proves $(1)\Rightarrow(2)$.

Conversely, suppose that $w_T^{\mathrm{loc}}(\eps)=O(\eps^\beta)$ as
$\eps\downarrow0$. Lemma~\ref{lem:entropy-from-local-width} gives
\[
\log N(T,\eps)
\leqslant C\int_\eps^\infty
\frac{w_T^{\mathrm{loc}}(t)^2}{t^3}\,dt
=O\bigl(\eps^{-2(1-\beta)}\bigr),
\]
and hence
\[
\mathcal D(T,u)
=O\left(\int_0^u\eps^{\beta-1}\,d\eps\right)
=O(u^\beta).
\]
Theorem~\ref{thm:dudley} gives a continuous modification $\widetilde X$ satisfying
\[
\limsup_{u\downarrow0}
\frac{\theta(\widetilde X,u)}{u^\beta}<\infty
\]
almost surely. As in the proof of Corollary~\ref{cor:holder}, this local
H\"older bound extends to all pairs of points in $T$, proving $(2)\Rightarrow(1)$.
\end{proof}

\section{Examples}\label{app:examples}

\begin{example}[Finite-dimensional sets]\label{ex:finite-dimensional}
Let $T$ be a compact set in a finite-dimensional subspace of
$H$. By monotonicity, $W(rT)\leqslant W(r(\conv T))$, and $\conv T$ is a
finite-dimensional compact convex set, so $W(r(\conv T))$ is a polynomial in
$r$. Hence $\rho(T)=0$, and Theorem~A gives $\eta(T)=0$ and
$\alpha_H(T)=1$. Concretely, if $G$ is the Gaussian vector representing the
process, then $X_t=\la G,t\ra$, and every sample path is Lipschitz on $T$ with
Lipschitz constant $\|G\|$.
\end{example}

\begin{example}[The Wiener spiral]\label{ex:wiener-spiral}
Let
\[
S=\left\{\mathbf 1_{[0,t]}:0\leqslant t\leqslant1\right\}\subset L^2[0,1].
\]
Since
\[
\la\mathbf 1_{[0,s]},\mathbf 1_{[0,t]}\ra=\min(s,t),
\qquad
\left\|\mathbf 1_{[0,t]}-\mathbf 1_{[0,s]}\right\|_{L^2}=|t-s|^{1/2},
\]
under the parametrization $t\mapsto\mathbf 1_{[0,t]}$, the isonormal process
on $S$ is a standard Brownian motion $(B_t)_{0\leqslant t\leqslant1}$. Thus a
modulus at Hilbert scale $\delta$ corresponds to one at time scale
$h=\delta^2$.

We compute the Wills functional. Put $\mu:=r\sqrt{\pi/2}$ and
\[
M_\mu:=\sup_{0\leqslant t\leqslant1}(B_t-\mu t).
\]
Then
\[
W(rS)=\mathbb E e^{2\mu M_\mu}.
\]
Indeed, since $\|\mathbf 1_{[0,t]}\|_{L^2}^2=t$, we have
\[
\sqrt{2\pi}\,rB_t-\pi r^2t=2\mu(B_t-\mu t),
\]
and the claim follows from the finite-subset definition
\eqref{eq:wills-gaussian} by continuity of Brownian paths.
% Note that the
% expression under the supremum is the logarithm of the exponential martingale
% of $B$.
The Bachelier--L\'evy formula \cite{RS70} yields, for $a\geqslant0$,
\[
\mathbb P(M_\mu\geqslant a)
=\overline\Phi(a+\mu)+e^{-2\mu a}\Phi(\mu-a),
\]
where $\Phi$ is the standard normal distribution function,
$\overline\Phi=1-\Phi$, and $\varphi=\Phi'$. Since $M_\mu\geqslant0$, the
tail-integral identity gives
\[
W(rS)
=1+2\mu\int_0^\infty
e^{2\mu a}\mathbb P(M_\mu\geqslant a)\,da.
\]

Integration by parts gives
\begin{align*}
2\mu\int_0^\infty e^{2\mu a}\overline\Phi(a+\mu)\,da
&=-\overline\Phi(\mu)
  +\int_0^\infty e^{2\mu a}\varphi(a+\mu)\,da\\
&=-\overline\Phi(\mu)+\int_0^\infty\varphi(a-\mu)\,da\\
&=2\Phi(\mu)-1,
\end{align*}
where we used $e^{2\mu a}\varphi(a+\mu)=\varphi(a-\mu)$ and
$\lim_{a\to\infty}e^{2\mu a}\overline\Phi(a+\mu)=0$.
For the second term, the substitution $x=\mu-a$ and the identity
$\int_{-\infty}^\mu\Phi(x)\,dx=\mu\Phi(\mu)+\varphi(\mu)$ yield
\begin{align*}
2\mu\int_0^\infty\Phi(\mu-a)\,da
&=2\mu\int_{-\infty}^\mu\Phi(x)\,dx\\
&=2\mu\bigl(\mu\Phi(\mu)+\varphi(\mu)\bigr).
\end{align*}
Combining the two terms gives
\[
W(rS)=2(1+\mu^2)\Phi(\mu)+2\mu\varphi(\mu).
\]
Substituting $\mu=r\sqrt{\pi/2}$, we obtain
\begin{equation}\label{eq:wiener-spiral-wills}
W(rS)
=\left(2+\pi r^2\right)\Phi\left(r\sqrt{\frac{\pi}{2}}\right)
+re^{-\pi r^2/4}.
\end{equation}
Since $\Phi\left(r\sqrt{\pi/2}\right)\to1$ and $re^{-\pi r^2/4}\to0$ as
$r\to\infty$, we conclude that $W(rS)=\pi r^2\left(1+o(1)\right)$. Thus the
Wills functional of the spiral grows polynomially, and
\[
\rho(S)=0.
\]
Theorem~A therefore gives $\eta(S)=0$ and $\alpha_H(S)=1$. Since
$\|\mathbf 1_{[0,t]}-\mathbf 1_{[0,s]}\|_{L^2}=|t-s|^{1/2}$, the exponent
$\alpha_H(S)=1$ with respect to the Hilbert metric corresponds to the classical
critical H\"older exponent $1/2$ of Brownian motion in the time parameter~$t$.
However, these exponent identities alone do not determine the logarithmic
correction in the exact modulus of continuity.

The full Wills estimate contains more information. Since
$\log W(rS)\leqslant C_0\log r$ for all $r\geqslant2$, choosing
$r_\eps=\eps^{-1}\sqrt{\log(1/\eps)}$ in
Theorem~\ref{thm:entropy-from-wills} gives, for all sufficiently small $\eps$,
\[
\sqrt{\log N(S,\eps)}
\leqslant
C\left(\frac{C_0\log r_\eps}{r_\eps\eps}+r_\eps\eps\right)
\leqslant
C_1\sqrt{\log(1/\eps)}.
\]
Hence
\[
\log N(S,\eps)\leqslant C_2\log(1/\eps),
\qquad \eps\to 0.
\]
This is the sharp entropy order. Indeed, a Hilbert ball of radius $\eps$
corresponds to a time interval of radius $\eps^2$, and a direct covering
argument gives
\[
N(S,\eps)=\left\lceil\frac{1}{2\eps^2}\right\rceil,
\qquad
\log N(S,\eps)\sim2\log\frac1\eps.
\]
Consequently,
\[
\mathcal{D}(S,\delta)
=
\int_0^\delta\sqrt{\log N(S,\eps)}\,d\eps
\leqslant
C_3\delta\sqrt{\log(1/\delta)},
\]
and Dudley's theorem therefore gives
\[
\theta(X,\delta)
=O\left(\delta\sqrt{\log(1/\delta)}\right)
\]
almost surely. In the time parameter this becomes
\[
\sup_{|t-s|\leqslant h}|B_t-B_s|
=O\left(\sqrt{h\log(1/h)}\right)
\qquad\text{almost surely as }h\downarrow0.
\]
Finally, L\'evy's exact modulus \cite[p.~172]{Levy54} states that
\[
\lim_{h\downarrow0}
\frac{\sup_{|t-s|\leqslant h}|B_t-B_s|}{\sqrt{2h\log(1/h)}}=1
\]
almost surely. Thus this order is sharp. The exponents $\rho(S)$, $\eta(S)$,
and $\alpha_H(S)$ are insensitive to logarithmic factors, whereas the
quantitative bound of Theorem~\ref{thm:entropy-from-wills} recovers
the sharp order of the metric entropy and, via Dudley's theorem, the sharp
order of the uniform modulus of continuity.
\end{example}

\begin{example}[The Brownian motion body]\label{ex:brownian-motion-body}
Let $S$ be the Wiener spiral from Example~\ref{ex:wiener-spiral} and set
$K:=\overline{\conv}S$, the closed convex hull of $S$. Gao and
Vitale~\cite{GV01} computed its intrinsic
volumes:
\[
V_j(K)=\frac{\kappa_j}{j!},
\qquad j\geqslant1.
\]
Since $\kappa_j=\pi^{j/2}/\Gamma(j/2+1)$, Stirling's formula gives
$\log(1/\kappa_j)=\frac12 j\log j+O(j)$ and $\log(j!)=j\log j+O(j)$, whence
\[
\log\frac1{V_j(K)}
=\frac32 j\log j+O(j).
\]
The coefficient formula \eqref{eq:order-from-coefficients} therefore yields
\[
\rho(K)=\frac23
\]
(and the accompanying type formula in Appendix~\ref{app:entire} shows that
$f_K$ has finite positive type $\sigma=\frac32(2\pi)^{1/3}$, so
$\log W(rK)\asymp r^{2/3}$).
This gives an example of intermediate order $0<\rho(K)<1$, for which Theorem~A
yields the exponents
\[
\eta(K)=1,
\qquad
\alpha_H(K)=\frac12.
\]
Recall that in $L^2[0,1]$ the closed convex hull $K=\overline{\conv}S$ coincides
with the class of all nonincreasing functions $f\colon[0,1]\to[0,1]$. Thus the
value $\eta(K)=1$ obtained here from the intrinsic volumes via Theorem~A agrees
with the classical result on the $L^2$-metric entropy of bounded monotone
functions, $\log N(K,\eps)\asymp\eps^{-1}$ as $\eps\to0$; see, for
example,~\cite[Theorem~5.2]{BS67} and~\cite[Theorem~2.7.5]{vdVW96}.
\end{example}

\section{Order and type of entire functions}\label{app:entire}

Let $F(z)=\sum_{j\geqslant0}a_jz^j\not\equiv0$ be an entire function and
$M_F(r)=\max_{|z|=r}|F(z)|$. If $F$ is nonconstant, its order and, when
$0<\rho<\infty$, its type are
\[
\rho=\limsup_{r\to\infty}\frac{\log\log M_F(r)}{\log r},
\qquad
\sigma=\limsup_{r\to\infty}\frac{\log M_F(r)}{r^\rho}.
\]
A nonzero constant is assigned order zero, and every polynomial has order zero.
For an entire function the order can be recovered from the
Taylor coefficients by
\begin{equation}\label{eq:order-from-coefficients}
\rho=\limsup_{\substack{j\to\infty\\a_j\neq0}}
\frac{j\log j}{\log(1/|a_j|)},
\end{equation}
and, if $0<\rho<\infty$, the type is given by
\[e\rho\sigma=\limsup_{j\to\infty} j |a_j|^{\rho/j}.\]
See, for instance, \cite[Ch.~1]{Lev96}.

%% file: Wills5_5.bbl
\begin{thebibliography}{99}

\bibitem{AHY21} D. Alonso-Guti\'errez, M. A. Hern\'andez Cifre, and J. Yepes Nicol\'as, \emph{Further inequalities for the (generalized) Wills functional}, Commun. Contemp. Math. \textbf{23} (2021), no. 3, 2050011.

\bibitem{AMS04} S. Artstein, V. D. Milman, and S. J. Szarek, \emph{Duality of metric entropy}, Ann. of Math. (2) \textbf{159} (2004), no. 3, 1313--1328, DOI 10.4007/annals.2004.159.1313.

\bibitem{BS67} M. S. Birman and M. Z. Solomyak, \emph{Piecewise-polynomial approximations of functions of the classes $W_p^\alpha$}, Mat. Sb. (N.S.) \textbf{73(115)} (1967), no. 3, 331--355; English transl. in Math. USSR-Sb. \textbf{2} (1967), no. 3, 295--317.

\bibitem{Bog98} V. I. Bogachev, \emph{Gaussian Measures}, Math. Surveys Monogr., vol. 62, Amer. Math. Soc., Providence, RI, 1998.

\bibitem{BZ88} Yu. D. Burago and V. A. Zalgaller, \emph{Geometric Inequalities}, Grundlehren Math. Wiss., vol. 285, Springer, Berlin, 1988; Russian original: Nauka, Leningrad, 1980.

\bibitem{Chev76} S. Chevet, \emph{Processus Gaussiens et volumes mixtes}, Z. Wahrscheinlichkeitstheorie verw. Gebiete \textbf{36} (1976), no. 1, 47--65.

\bibitem{DGZ25} M. Dospolova, M. Germanskov, and D. Zaporozhets, \emph{On Steiner entire function}, arXiv:2507.11626 (2025).

\bibitem{Dud67} R. M. Dudley, \emph{The sizes of compact subsets of Hilbert space and continuity of Gaussian processes}, J. Funct. Anal. \textbf{1} (1967), 290--330.

\bibitem{Dud73} R. M. Dudley, \emph{Sample functions of the Gaussian process}, Ann. Probab. \textbf{1} (1973), no. 1, 66--103.

\bibitem{Fer75} X. Fernique, \emph{R\'egularit\'e des trajectoires des fonctions al\'eatoires gaussiennes}, \'Ecole d'\'Et\'e de Probabilit\'es de Saint-Flour IV---1974, Lecture Notes in Math., vol. 480, Springer, Berlin, 1975, pp. 1--96.

\bibitem{GV01} F. Gao and R. A. Vitale, \emph{Intrinsic volumes of the Brownian motion body}, Discrete Comput. Geom. \textbf{26} (2001), no. 1, 41--50.

\bibitem{Had75} H. Hadwiger, \emph{Das Wills'sche Funktional}, Monatsh. Math. \textbf{79} (1975), 213--221.

\bibitem{HHC08} M. Henk and M. A. Hern\'andez Cifre, \emph{Notes on the roots of Steiner polynomials}, Rev. Mat. Iberoam. \textbf{24} (2008), no. 2, 631--644.

\bibitem{HCYN13} M. A. Hern\'andez Cifre and J. Yepes Nicol\'as, \emph{On the roots of the Wills functional}, J. Math. Anal. Appl. \textbf{401} (2013), no. 2, 733--742.

\bibitem{IN68} K. It\^o and M. Nisio, \emph{On the oscillation functions of Gaussian processes}, Math. Scand. \textbf{22} (1968), 209--223.

\bibitem{KL93} J. Kuelbs and W. V. Li, \emph{Metric entropy and the small ball problem for Gaussian measures}, J. Funct. Anal. \textbf{116} (1993), no. 1, 133--157, DOI 10.1006/jfan.1993.1107.

\bibitem{Lev96} B. Ya. Levin, \emph{Lectures on Entire Functions}, Transl. Math. Monogr., vol. 150, Amer. Math. Soc., Providence, RI, 1996.

\bibitem{Levy54} P. L\'evy, \emph{Th\'eorie de l'addition des variables al\'eatoires}, 2nd ed., Gauthier-Villars, Paris, 1954.

\bibitem{Lif95} M. A. Lifshits, \emph{Gaussian Random Functions}, Math. Appl., vol.~322, Kluwer, Dordrecht, 1995.

\bibitem{LL99} W. V. Li and W. Linde, \emph{Approximation, metric entropy and small ball estimates for Gaussian measures}, Ann. Probab. \textbf{27} (1999), no. 3, 1556--1578, DOI 10.1214/aop/1022677459.

\bibitem{Lif12} M. A. Lifshits, \emph{Lectures on Gaussian Processes}, SpringerBriefs in Mathematics, Springer, Heidelberg, 2012, DOI 10.1007/978-3-642-24939-6.

\bibitem{MR06} M. B. Marcus and J. Rosen, \emph{Markov Processes, Gaussian Processes, and Local Times}, Cambridge Stud. Adv. Math., vol. 100, Cambridge Univ. Press, Cambridge, 2006.

\bibitem{McM91} P. McMullen, \emph{Inequalities between intrinsic volumes}, Monatsh. Math. \textbf{111} (1991), no. 1, 47--53.

\bibitem{MP87} V. D. Milman and G. Pisier, \emph{Gaussian processes and mixed volumes}, Ann. Probab. \textbf{15} (1987), no. 1, 292--304.

\bibitem{Mou25} J. Mourtada, \emph{Universal coding, intrinsic volumes, and metric complexity}, J. Eur. Math. Soc. (2025), published online first, DOI 10.4171/JEMS/1649; arXiv:2303.07279.

\bibitem{Mou26} J. Mourtada, \emph{Quantitative characterization of Gaussian small ball probabilities}, arXiv:2611.17071 (2026).

\bibitem{RS70} H. Robbins and D. Siegmund, \emph{Boundary crossing probabilities for the Wiener process and sample sums}, Ann. Math. Statist. \textbf{41} (1970), 1410--1429.

\bibitem{Sch14} R. Schneider, \emph{Convex Bodies: The Brunn--Minkowski Theory}, 2nd ed., Encyclopedia Math. Appl., vol. 151, Cambridge Univ. Press, Cambridge, 2014.

\bibitem{Sud71} V. N. Sudakov, \emph{Gaussian random processes and measures of solid angles in Hilbert space}, Dokl. Akad. Nauk SSSR \textbf{197} (1971), 43--45; English transl. in Soviet Math. Dokl. \textbf{12} (1971), 412--415.

\bibitem{Sud76} V. N. Sudakov, \emph{Geometric problems in the theory of infinite-dimensional probability distributions}, Trudy Mat. Inst. Steklov. \textbf{141} (1976); English transl. in Proc. Steklov Inst. Math. \textbf{141} (1979), 1--178.

\bibitem{Tal87} M. Talagrand, \emph{Regularity of Gaussian processes}, Acta Math. \textbf{159} (1987), 99--149.

\bibitem{Tal21} M. Talagrand, \emph{Upper and Lower Bounds for Stochastic Processes: Decomposition Theorems}, 2nd ed., Ergeb. Math. Grenzgeb. (3), vol. 60, Springer, Cham, 2021.

\bibitem{Tsi85} B. S. Tsirelson, \emph{A geometric approach to maximum likelihood estimation for an infinite-dimensional Gaussian location. II}, Teor. Veroyatnost. i Primenen. \textbf{30} (1985), no. 4, 772--779; English transl. in Theory Probab. Appl. \textbf{30} (1986), 820--828.

\bibitem{vdVW96} A. W. van der Vaart and J. A. Wellner, \emph{Weak Convergence and Empirical Processes: With Applications to Statistics}, Springer Ser. Statist., Springer, New York, 1996.

\bibitem{Vit96} R. A. Vitale, \emph{The Wills functional and Gaussian processes}, Ann. Probab. \textbf{24} (1996), no. 4, 2172--2178.

\bibitem{Vit01} R. A. Vitale, \emph{Intrinsic volumes and Gaussian processes}, Adv. Appl. Probab. \textbf{33} (2001), no. 2, 354--364.

\bibitem{Vit09} R. A. Vitale, \emph{On the oscillation of a Gaussian process}, High Dimensional Probability V: The Luminy Volume, Inst. Math. Stat. Collect., vol. 5, IMS, Beachwood, OH, 2009, pp. 219--221.

\bibitem{Wil73} J. M. Wills, \emph{Zur Gitterpunktanzahl konvexer Mengen}, Elem. Math. \textbf{28} (1973), 57--63.

\end{thebibliography}
